\documentstyle[12pt]{article}

\begin{document}

\newtheorem{theorem}{Theorem}[section]
\newtheorem{corollary}[theorem]{Corollary}
\newtheorem{definition}[theorem]{Definition}
\newtheorem{conjecture}[theorem]{Conjecture}
\newtheorem{question}[theorem]{Question}
\newtheorem{lemma}[theorem]{Lemma}
\newtheorem{proposition}[theorem]{Proposition}
\newtheorem{example}[theorem]{Example}
\newtheorem{problem}[theorem]{Problem}
\newenvironment{proof}{\noindent {\bf
Proof.}}{\rule{3mm}{3mm}\par\medskip}
\newcommand{\remark}{\medskip\par\noindent {\bf Remark.~~}}
\newcommand{\pp}{{\it p.}}
\newcommand{\de}{\em}

\newcommand{\JEC}{{\it Europ. J. Combinatorics},  }
\newcommand{\JCTB}{{\it J. Combin. Theory Ser. B.}, }
\newcommand{\JCT}{{\it J. Combin. Theory}, }
\newcommand{\JGT}{{\it J. Graph Theory}, }
\newcommand{\ComHung}{{\it Combinatorica}, }
\newcommand{\DM}{{\it Discrete Math.}, }
\newcommand{\ARS}{{\it Ars Combin.}, }
\newcommand{\SIAMDM}{{\it SIAM J. Discrete Math.}, }
\newcommand{\SIAMADM}{{\it SIAM J. Algebraic Discrete Methods}, }
\newcommand{\SIAMC}{{\it SIAM J. Comput.}, }
\newcommand{\ConAMS}{{\it Contemp. Math. AMS}, }
\newcommand{\TransAMS}{{\it Trans. Amer. Math. Soc.}, }
\newcommand{\AnDM}{{\it Ann. Discrete Math.}, }
\newcommand{\NBS}{{\it J. Res. Nat. Bur. Standards} {\rm B}, }
\newcommand{\ConNum}{{\it Congr. Numer.}, }
\newcommand{\CJM}{{\it Canad. J. Math.}, }
\newcommand{\JLMS}{{\it J. London Math. Soc.}, }
\newcommand{\PLMS}{{\it Proc. London Math. Soc.}, }
\newcommand{\PAMS}{{\it Proc. Amer. Math. Soc.}, }
\newcommand{\JCMCC}{{\it J. Combin. Math. Combin. Comput.}, }
\newcommand{\GC}{{\it Graphs Combin.}, }

\title{The  Laplacian Spectra of  Graphs and Complex Networks\thanks{
This work is supported by National Natural Science Foundation of
China (No:10971137), the National Basic Research Program (973) of
China (No.2006CB805900)  and a grant of Science and Technology
Commission of Shanghai Municipality (STCSM, No: 09XD1402500).\vskip
0.1cm 
 \vskip 0.1cm
Corresponding author: Xiao-Dong Zhang (xiaodong@sjtu.edu.cn)
 }}
\author{ Ya-Hong Chen\\
{\small College of Education, Lishui University, Lishui, Zhejiang
323000, P.R. China}\\
Rong-Ying Pan\\
{\small Suzhou Vocational University Suzhou, Jiangsu, 210000, P.R.
China}\\
Xiao-Dong Zhang \\
{\small Department of Mathematics, Shanghai Jiao Tong University} \\
{\small  800 Dongchuan road, Shanghai, 200240, P.R. China}\\
{\small Email:  xiaodong@sjtu.edu.cn}\\
 {\small \em  Dedicated to
Professor Jiong-Sheng Li on the occasion of his 75th birthday}
  }
\maketitle
 \begin{abstract}
The paper is a  brief survey of some recent new results and progress
of the  Laplacian spectra of graphs and complex networks (in
particular, random graph and the small world network). The main
contents contain the spectral radius of the graph Laplacian for
given a degree sequence, the Laplacian coefficients, the algebraic
connectivity and the graph doubly stochastic matrix, and the spectra
of random graphs and the small world networks. In addition, some
questions are proposed.
 \end{abstract}


\vskip 0.5cm

\section{Introduction}

 The study of graph Laplacian spectrum realized increasingly
 connections with many other areas.  The objects arise in very diverse application, from
 combinatorial optimization to  differential geometry, mathematical biology,
 computer science, machine learning, etc.)
 The Laplacian spectrum can be used to extract useful and important information
 about some graph invariants (for examples, expansion and vertex partition)
 that are hard to computer or estimate (see \cite{alon1985}).  One of the most
 fascinating aspects of  applications of eigenvalue methods in Combinatorics
   is the spectrum appear as
 a tool to prove results that appear to have nothing to do with the
 spectrum itself (see the excellent monographs \cite{godsil2001}).
 The smallest nonzero and largest eigenvalues \cite{mohar1997} can
 be  expressed as solutions to a quadratic optimization problem. An
 important use of eigenvalues is Lov\'{a}sz's notation of the
 $\vartheta$ function (\cite{lovasz1975}) which was initiated by problems
  in communication networks.

    The graph Laplacian can be used in a number of ways to provide
    interesting geometric representations of a graph. One of the
    most work (\cite{colin1990, godsil2001}) is the Colin de Verdi\`{e}re number of a graph, which
    is regarded as one of the most important recent developments in
    graph theory. The parameter, which is minor-monotone, determines
    the embeddability properties of graphs.

Alon and Milman \cite{alon1985} proved that the Laplacian spectrum
of graphs plays a crucial role  from the explicit constructions of
expander graphs and superconcentrators to the design of various
randomized algorithms.  Based on the Laplacian eigenvalues and
isoperimetric properties of graphs, Lubotzky, Phillips and Sarnak in
\cite{lubotzky1988} gave several family explicit constructions of
good expander graphs.

The Laplacian matrix of graphs can be regarded as  the discrete
Laplacian operator on the differential manifolds. There is an
interesting bilateral link between spectral graph theory and
spectral Riemannian geometry. The concepts and methods of spectral
geometry make it possible to obtain new results in the study of
graph Laplacian spectrum. Conversely, the results in spectral graph
theory may be transferred  to the Laplacian operator on manifolds.
We refer to the excellent book \cite{Chung1998}.

 At the dawn of the new century, the power law networks (the scale-free networks)
  \cite{barabasi1999} and the  small world networks
  \cite{watts1998}
were discovered and studied. Since then, the analysis and modeling
of networks, and networked dynamical system, have been the subject
of considerable interdisciplinary interests, including physics,
mathematics, computer science, biology, economics.  The results have
been  called the "new science of networks" (see \cite{watts2004}).

The Laplacian spectra of  a various way of modeling of networks
which represent the real networks can be used to give a tentative
classification scheme for empirical networks  and  provide useful
insight into the relevant structural properties of real networks.
Spectral techniques and methods based on the analysis of the largest
eigenvalues and eigenvectors of some complex networks (for example
the web networks \cite{newman2003}) have proven algorithmically
successful in detecting  communities and clusters.
 Mihail and Papadimitriou \cite{mihail2004} showed that the largest
eigenvalue of a power law graph with exponent $\beta$ has power law
distribution if the exponent $\beta$ of the power law graph
satisfies $\beta>3$. Chung  etc. in \cite{chung 2003,chung 2003b}
established relationships between the spectrum of graphs and the
structure of complex networks.  It is showed that spectral methods
are central in detecting clusters and finding patterns in various
applications. For more information, we refer to  the monograph
\cite{chung 2006} and a survey \cite{newman2003}.

In this paper, we survey recent progress of  the Laplacian spectra
of graphs and complex networks.  For basic information and earlier
results, we refer several books and  survey papers  such as
\cite{Chung1998, chung 2006, merris1994, merris1995, mohar1991,
mohar1997, zhang2007}, for a detailed introduction and applications.
The rest of this paper is organized as follows. In Section 2, some
notations and results are presented. In Section 3, the Laplacian
spectral radius of graphs with given degree sequences is extensively
investigated. In Section 4, the relationship between the Laplacian
coefficient and the ordering of graphs are studied, in particular,
the Mohar's problems \cite{mohar2008} on the topic are discussed. In
Section 5, we deeply studied the algebraic connectivity and
structure and properties of the graph doubly stochastic matrix, in
particular, on Merris's problems \cite{Merris1998b} In Section 5,
the Laplacian spectrum of random graphs and small world networks are
discussed.

\section{Preliminary}

Let $G = (V(G),~E(G))$ be a simple graph (no loops or multiple
edges) with vertex set $V(G)=\{v_1,\cdots, v_n\}$ and edge set
$E(G)$. The {\it degree} of vertex $v_i$, denote by $d(v_i)$ or
$d_G(v_i)$, is the number of the edges incident with $v_i$. Let
$D(G)=diag(d(u), u\in V)$ be the diagonal matrix of vertex degrees
of $G$ and $A(G)=(a_{ij})$ be the $(0,1)$ {\it  adjacency matrix} of
$G$, where $a_{ij}=1$ for $v_i$ adjacent to $v_j$ and 0 for
elsewhere.  Then the matrix $L(G)=D(G)-A(G)$ is called the {\it
Laplacian matrix} of a graph $G$. It is obvious that $L(G)$ is
positive semidefinite and singular $M-$matrix. Thus the set of  all
eigenvalues of $L(G)$  are called the {\it Laplacian spectrum}  of
$G$ and arranged in nonincreasing order:
$$\lambda_1\ge \lambda_2\ge\cdots\ge \lambda_{n-1}\ge
\lambda_n=0.$$ When more than one graph is under discussion, we
may write $\lambda_i(G)$ instead of $\lambda_i$.
 From the
matrix-tree theorem, it is easy to see that $\lambda_{n-1}>0$ if and
only if $G$ is connected. This observation led M. Fiedler
\cite{fiedler1973} to define the algebraic connectivity of $G$ by
$\alpha(G)=\lambda_{n-1}(G)$, which can  be  used as a quantitative
measure of connectivity. Further, the characteristic polynomial of
$L(G)$ is called the {\it Laplacian polynomial} of $G$ and is
denoted by
$$\Phi(G)=\det(\lambda
I_n-L(G))=\sum_{i=0}^n(-1)^ic_i\lambda^{n-i}.$$ It is easy to see
that $c_0=1, $ $ c_1=2|E(G)|$, $ c_{n}=0$ and $c_{n-1}=n\tau(G)$,
where $\tau(G)$ is the number of spanning trees of $G$. It is well
known that
$$\sum\lambda_{i_1}\cdots\lambda_{i_k}=c_k,$$
where the sum is taken over all subsets of $k$ elements in the set
$\{1, \cdots, n\}$. Hence the spectrum $ \{\lambda_1, \cdots,
\lambda_n\}$ and the set consisting of $c_0, \cdots, c_n$ are
mutually determined.

\section{The Laplacian spectral radius and graphic degree sequences}
 A nonincreasing sequence of nonnegative integers
 $\pi=(d_0,d_1,\cdots, d_{n-1})$ is called {\it graphic} if there
 exists a simple graph with $\pi$ as its  vertex degree sequence.
Moreover, $\pi$ is called a unicyclic degree sequence, if  there
exists a unicyclic graph (i.e., the connected graph with only one
cycle) with $\pi$.  From the point of view of linear algebra, many
properties of the graph
 Laplacian spectrum can be characterized by  simple graph
 invariants, such as vertex degree. There are a lot of results on
 the upper and lower bounds for the spectral radius of graphs in terms of
 the maximum degree, the minimum degree, etc.(see \cite{zhang2007} for the detail).
 Motivated by the recent results in terms of vertex degrees,  Zhang in \cite{zhang2008}
  generally
 propose the following question.
\begin{problem}\cite{zhang2008} For a given  graphic degree sequence $\pi$, let
\begin{eqnarray*}
 {\mathcal G}_{ \pi}=\{ G \ | \   G \ {\rm is\  connected \ with\  \pi \ as
 \  its \  degree\  sequence}\}.
 \end{eqnarray*}
 Study the  sharp upper (lower) bounds  for the Laplacian spectral  radius of
all graphs $G$ in ${\mathcal G}_{ \pi}$ and characterize all
extremal graphs which attain the upper (lower) bounds.
\end{problem}
 It is natural to consider the tree degree sequence, since tree is a
 simplest connected graph.  With aid of the properties of the
 eigenvector corresponding to the spectral radius,
 \cite{zhang2008} gave the extremal tree in the set of given tree degree sequence which
 has the maximum Laplacian spectral radius.

Let $\pi=(d_0, d_1,\cdots, d_{n-1}) $  with $n\ge 3$ be  a given
nonincreasing degree sequence  of  some  tree. Now we construct a
special tree $T_{\pi}^*$ with degree sequence $\pi$ by using a
"breadth-first" scheme, which considers the $n$ vertices of $T^*$ to
be partitioned into a sequence of layers starting with the $0^{th}$
layer consisting of a single vertex $v_{01}$ of maximum degree. Then
recursively develop the layers, with the $1^{th}$ layer having $d_0$
vertices each connected to $v_{01}$. Denote the number of vertices
in the $m^{th}$ layer by $l_m$, and the number in all layers
preceding layer $m$ by $l_{<m}$. Given $l_m$ vertices in a layer
$m\ge 1,$ choose the degree of the $i^{th}$ vertices $v_{m, i}$ in
this layer to be $d_{i+l_{<m}}$ with all but 1 of its connections to
be to vertices $v_{m+1, k}$ in the next layer for
$k=j+\sum_{p=1}^{i-1}(d_{p+l_{<m}}-1), $ with $j=1$ to
$d_{i+l_{<m}}-1$. The number of vertices in this next layer then is
$l_{m+1}=\sum_{p=1}^{l_m}(d_{p+l_{<m}}-1)$, and the process is
continued till all $n$ vertices are accounted for.

\begin{theorem}\cite{zhang2008}\label{max-degre}
For a given degree sequence $\pi$ of some tree, let
\begin{eqnarray*}
{\mathcal T}_{\pi}=\{ T \ | \   T \ {\rm is\  tree \ with\  \pi \ as
 \  its \  degree\  sequence}\}.
 \end{eqnarray*}
   Then $T_{\pi}^*$  is a unique
tree with the maximum Laplacian spectral radius in ${\mathcal T}_{
\pi}$.
\end{theorem}
 Let $\pi=(d_0,\cdots,
d_{n-1})$ and
 $\pi^{\prime}=(d_0^{\prime}, \cdots, d_{n-1}^{\prime})$ be two
 nonincreasing sequences. If
 $\sum_{i=0}^{k}d_i\le\sum_{i=0}^kd_i^{\prime}$ and
 $\sum_{i=0}^{n-1}d_i=\sum_{i=0}^{n-1}d_i^{\prime}$, then the
 sequence $\pi^{\prime}$ is said to {\it major} the sequence $\pi$ and denoted by
 $\pi\triangleleft \pi^{\prime}$.  Further, the majorization theory on the two different tree
 degree sequences.
\begin{theorem}\cite{zhang2008}
Let $\pi$ and $\tau$ be two different tree degree sequences with the
same order. Let $T_{\pi}^*$ and $T_{\tau}^*$ have the maximum
 Laplacian spectral radii in ${\mathcal T}_{ \pi}$ and ${\mathcal
T}_{\tau}$, respectively.
 If
$\pi\triangleleft \tau$, then
$\lambda(T_{\pi}^*)<\lambda(T_{\tau}^*)$.
\end{theorem}

Next, it will be interesting to characterize the extremal graphs
with the set of all unicyclic graphs for given unicyclic graphic
degree sequence $\pi$. But it seems to be more difficult.
 However, there are some partial results on the unicyclic graphs.
 For examples,  Tan and Zhang in \cite{tan2003}  presented a sharp bound for the spectral radius of the Laplacian
 Matrix of unicyclic graphs of order $n$ and the matching number $\beta$. Further,
  they  characterized all extremal graphs which attained the upper
  bounds.
  \begin{theorem}\cite{tan2003}
  Let $G$ be any unicyclic graph of order $n$ with the matching
  number $\beta$. Then the Laplacian spectral  radius $\lambda_1(G)$
  of $G$  is no more than the the largest root of the following
  equation
  $$\lambda^3-(n-\beta+5)\lambda^2+(3n-3\beta+7)\lambda-n=0$$
 with equality if and only if $G$ is obtained from a triangle and
 $\beta-2$ paths with length 2 and $n-2\beta+1$ edges by identifying
 one end vertex of them.
 \end{theorem}

 \section{Graph Laplacian coefficients }

  It is an interesting and important problem how to order graphs for a given set of graphs,
   which has many application in computer science.
   There are many approaches to ordering these graphs. For example, the
graphs can be ordered
  lexicographically according to their eigenvalues in nonincreasing order
  (for example, see \cite{Cvetkovic1980}, pp.268-269 and \cite{cvetkovic1988}, p.70).
  Grone and Merris in \cite{grone1990a} used the algebraic connectivity of a tree $T$
   to    order trees. Recently,  Zhang in \cite{zhang2007b} further
investigated ordering trees by their algebraic connectivity.

On the other hand, Another graph invariant, the Wiener index,  can
also be used to order graphs. The {\it Wiener index} of a connected
graph is the sum of all distances between unordered pairs of
vertices of a connected graph. In other words,
$$W(G)=\sum_{\{v_i, v_j\}\subseteq V(G)}d_G(v_i, v_j),$$
where $d_G(v_i, v_j)$ is the  distance between vertices $v_i$ and
$v_j$, i.e.,  the minimum number of edges between $v_i$ and $v_j$.
 The
Wiener index was introduced in 1947 and  extensively studied by
chemists
 and mathematicians (for example, see \cite{dobrynin2001, yan2005} for more detail).
 Since it is correlated
with several graph properties and some physical and chemical
properties of molecular graphs,  the Wiener index can be  used to
order trees (see \cite{dong2006}) and extract useful information
\cite{zhou2008} from the structure of graphs.
 It is known that the Wiener index of a tree is equal to $c_{n-2}$ (for example, see
\cite{merris1989} or \cite{mohar1991}). Thus trees with the same
Wiener index may be further ordered by other Laplacian coefficients.

Recently, Mohar in \cite{mohar2007} and his homepage
\cite{mohar2008} proposed some
 problems on how to order trees with the Laplacian coefficients.
Let $T_1$ and $T_2$ be two trees of order $n$. Let  $ r$ (resp. $s$)
be the smallest (resp. largest) integer such that $c_r(T_1)\neq
c_r(T_2)$ (respectively, $c_s(T_1)\neq c_s(T_2)$).  Thus $r$ and $s$
exist if and only if $T_1$ and  $T_2$ are not Laplacian cospectral.
Hence two partial orderings can be defined as follows:  If
$c_r(T_1)<c_r(T_2), $     $T_1$ is ``{\it smaller than} " $T_2$ and
denote $T_1\prec^1 T_2$. If $c_s(T_1)<c_s(T_2), $  $T_1$ is ``{\it
smaller than} " $T_2$ and denote $T_1\prec^2 T_2$.
 Let ${\mathcal{T}}_n$ be the set of all trees of order $n$. Another partial ordering in
   ${\mathcal{T}}_n$  can be defined as follows: For
  any two trees  $T_1, T_2\in {\mathcal{T}}_n$, if $(c_0(T_1),\cdots, c_n(T_1))\le
   (c_0(T_2),\cdots, c_n(T_2))$, i.e.,  $c_i(T_1)\le
  c_i(T_2)$  with $i=0, \cdots, n$, we say that $T_1$ is {\it dominated} by
  $T_2$ and denote $T_1\preceq T_2$. Mohar in \cite{mohar2008} proposed the following
  questions:

\begin{problem}\cite{mohar2008}\label{prob1} Do there exist two trees $T_1$ and $T_2$ of order $n$
 such that
$ T_1\prec^1 T_2$ and $T_2\prec^2 T_1$?
\end{problem}

 \begin{problem}\cite{mohar2008}\label{prob2}
 Do there exist two trees $T_1$ and $T_2$ of order $n$ such that
 $T_1\prec^1 T_2$ and $T_1\prec^2 T_2$, but there is an index $i$
  such that $c_i(T_1) > c_i(T_2)?$
  \end{problem}

\begin{problem}\cite{mohar2008}\label{prob3}
Let ${\mathcal{T}}_n$ be the set of all trees of order $n$. How
large chains and anti-chains of pairwise non- Laplacian-cospectral
trees are there?
 \end{problem}

\begin{problem}\cite{mohar2008}\label{prob4} Let $T_1$ and $T_2$ be two trees of order $n$ with
$T_1 \preceq  T_2$. Let ${\mathcal{T}}(T_1, T_2)$ be the set of all
trees $T$ of order $n$ with $T_1\preceq  T \preceq T_2$. For which
trees $T_1$ and $T_2$ has ${\mathcal{T}}(T_1, T_2)$ only one minimal
element up to cospectrality, i.e., when are all minimal elements in
${\mathcal{T}}(T_1, T_2)$ cospectral?
\end{problem}

On the problems \ref{prob1} and \ref{prob2},  Zhang in
\cite{zhang2009b} gave a positive  answer for problems \ref{prob1}
and \ref{prob2} by presenting two examples, and  investigated all
majorization relationship among all trees of order $n$ with diameter
3. It may be more interesting to gave infinite family trees such
that the results holds and  to investigate structure  of trees with
the properties.

 Gutman and Pavlovic in \cite{gutman2003} proved there
exists only maximal and minimal elements in  the set of all trees
 of order $n$ with  respective to $({\mathcal T}_n, \preceq)$.
\begin{theorem}\cite{gutman2003}
Let $T$, $K_{1, n-1}$ and $P_n$ be any tree,  the star and path of
order $n$. Then
$$c_i(K_{1,n-1})\le c_i(T)\le c_i(P_n).$$
In other words, $K_{1, n-1} \preceq T \preceq P_n$.
\end{theorem}
In addition, Ili\'{c} in \cite{ilic2009} studied the maximal element
for given a set of all trees of order $n$ with fixed the maximum
degree.

\begin{theorem}\cite{ilic2009}
Let $T$ be any tree of order $n$ with the maximum degree $\Delta$.
Then $$ c_i(T)\le c_i(B_{n,\Delta}), \ \ i=0, \cdots, n$$
 where $B_{n, \Delta}$  is the tree
consisting of a star $K_{1, \Delta}$ and a path of length
$n-\Delta-1$ attached to an arbitrary pendent vertex of the star.
Moreover, equality holds for some $2\le i\le n-2$ if and only if $T$
is $B_{n, \Delta}$
\end{theorem}

 Recently, Wang \cite{wang2008} and Zhang
etc. \cite{zhang2008b} independently proved the minimum Wiener index
in the set of all trees with give a tree degree sequence.

\begin{theorem}\cite{wang2008, zhang2008b}\label{wiener}
For a given degree sequence $\pi$ of some tree, let
$${\mathcal T}_{\pi}=\{ T \ | \   T \ {\rm is\ a\  tree \ with\  \pi \ as
 \  its \  degree\  sequence}\}.$$  Then $T_{\pi}^*$ is a unique
tree with the minimum Wiener index in ${\mathcal T}({ \pi})$.
\end{theorem}
Motivated  Theorems~\ref{max-degre}, \ref{wiener} and the related
results, we proposed the following conjecture

\begin{conjecture}\label{con-coeff}
For a given degree sequence $\pi$ of some tree, let
$${\mathcal T}_{\pi}=\{ T \ | \   T \ {\rm is\ a\  tree \ with\  \pi \ as
 \  its \  degree\  sequence}\}.$$  Then for any $ T\in {\mathcal
 T}_{\pi}$, we have
 $$c_i(T_{\pi}^*)\le c_i(T), \ \ i=0, \cdots, n $$
 with equality holding for some $2\le i\le n-2$ if and only if  $T$ is $T_{\pi}^*$.
\end{conjecture}

{\bf Remark:} If Conjecture~\ref{con-coeff} holds, then
Theorems~\ref{max-degre} and \ref{wiener} will be direct corollary
of the conjecture.

In addition, Stevanovi\'{c} and Ili\'{c} \cite{stevanovic2009}
investigated the properties of the Laplacian coefficients of
unicyclic graphs. Moreover, The related results on the Laplacian
coefficients can be seen \cite{ilic2009b}.

\section{The algebraic connectivity and doubly stochastic matrix}
In the  study of chemical information processing,  Golender et al.
\cite{Goiender1981} introduced another important matrix: {\it doubly
stochastic graph matrix} which is related to the graph Laplacian
matrix and can  be used to describe properties of topological
structure of chemical molecular  graphs. Let $I_n$ be the $n\times
n$ identity matrix and $\Omega(G)=(I_n+L(G))^{-1}=(\omega_{ij})$. It
is easy to see (\cite{Goiender1981} or  \cite{Merris1997}) that
$\Omega(G)$ is a doubly stochastic matrix. Thus  $\Omega(G)$ is
called the {\it doubly stochastic graph matrix}. On the other hand,
Chebotarev in \cite{Chebotrarev1997}  pointed out that the doubly
stochastic graph matrix may be used to measure the proximity among
vertices and evaluate the group cohesion in the construction of
sociometric indices and represent a random walk. Moreover, with the
 entries of the doubly stochastic matrix, many topological measurements (such as
 dissociation, solitariness, provinciality, and  etc) of small social groups
  on the basis of given relations on them  can be quantified.
 Merris  in
 \cite{Merris1998b} studied relationship
  between the the algebraic connectivity and the entries of the
  graph doubly stochastic matrix. In particular, he proposed two
  conjectures and two problems.  One of the conjectures and one of
  the problems are follows:

\begin{conjecture}\cite{Merris1998b}
\label{conj1} Let $G$ be a graph on $n$ vertices. Then
$$\alpha(G)\ge 2(n+1)\omega(G),$$
where $\omega(G)$ is the smallest
 entry of   $\Omega(G)=(\omega_{ij})$,
i.e., $\omega(G)=\min\{\omega_{ij}; 1\le i, j\le n\}$.
\end{conjecture}

\begin{question}\cite{Merris1998b}\label{merris}
 Let $G$ be a  simple graph on  vertex set $V=\{v_1,\cdots,
v_n\}$ with doubly stochastic graph matrix
$\Omega(G)=(\omega_{ij})$. Let $\rho(v_i,
v_j)=\omega_{ii}+\omega_{jj}-2\omega_{ij}$ and
$$\underline{r}(G)=\min\{\sum_{j=1}^n \rho(v_i, v_j),\ 1\le i\le
n\}.$$ Does $d_k > d_i$ for all $i\neq k$, imply
$r(k)=\underline{r}(G)$?
\end{question}

Zhang and Wu  in \cite{zhang2005a}  presented an example to
 illustrate that Merris' Conjecture~\ref{conj1}  does not
 hold  generally.  Recently, Zhang in \cite{zhang2009c}  carefully
 investigated when this conjecture still hold and gave many examples to illustrate
 that this conjecture does not holds. In particular,

 \begin{theorem}\cite{zhang2009c}\label{diaalge}
Let $T$ be a tree of order $n\ge 4$ with diameter $d$. If $d\ge
\frac{\lg 3+3\lg n}{\lg(3+\sqrt{5})-\lg 2}-1$, then $\alpha(T)\ge
2(n+1)\omega(T)$.
\end{theorem}

\begin{theorem}\cite{zhang2009c}\label{upper}
Let $T$ be a tree of order $n$ with $p$ non-pendant  vertices. Then
$$\alpha(T)\ge \frac{(n+p)\omega(T)}{1-(n+p)\omega(T)}$$ with
equality if and only if $T$ is the star graph $K_{1, n-1}$.
\end{theorem}

{\bf Remark} From the above results, we may see that
Conjecture\ref{conj1} holds for many trees, since the diameter of
any random trees is almost equal $O(\lg n)$. While
Conjecture\ref{conj1} does not holds for  smaller diameter and
larger order.

On Merris question~\ref{merris}, Zhang in \cite{zhang2010} proved
the following results

\begin{theorem}\cite{zhang2010}\label{main}
\hskip 20cm (i)  (i) There exists a family of graphs  with
$d(v_k)>d(v_i)$ for  all $i\neq k$  but $r(k)>\underline{r}(G)$.

\noindent (ii) Let $G$ be a simple graph. If  $d(v_k) \ge 2d(v_i)$
for all $i\neq k$, then  $r(k)=\underline{r}(G)$.
\end{theorem}
 {\bf Remark:}  There are  many examples to illustrate the question \ref{merris}
 still holds for $d(v_k)\ge 2d(v_i)-2$.  Further Zhang in \cite{zhang2010} proposed
 following question:
\begin{question}\cite{zhang2010}
\label{conj} Let $G$ be a  simple connected graph on $n$ vertices
$\{v_1,\dots, v_n\}$ with the doubly stochastic graph matrix
$\Omega(G)=(\omega_{ij})$. If $d(v_k)=\max \{ d(v_1), d(v_2),\cdots,
d(v_n)\}$ and $d(v_k)\ge 2d(v_i)-2$ for all $i\neq k$, does
$r(k)=\underline{r}(G)$ hold ?
\end{question}
The  related works on the algebraic connectivity and the structure
of the graph doubly stochastic matrix can be
 referred to \cite{berman2000, Chebotrarev2002,
zhang2005b, zhang2005a}.
\section{ The algebraic connectivity of random graphs}

Recently, there has been much interest in studying the small world
networks and attempting to model their properties using random
graphs. The study of complex systems in terms of random graphs was
initiated by Erd\"{o}s and R\'{e}nyi \cite{erdos1960b}. Recently
Watts and Strogatz \cite{watts1998} introduced the small-world
network by interpolating order and random.

 For given $0\le p\le 1$, a {\it ER random graph} ${\mathcal{G}}(n,p)$
  is a graph of $n$ nodes
and edges between nodes occurs independently with probability $p$.
This random network model has been intensively studied (see
\cite{bollobas2001} and the references therein). Juh\'{a}sz in
\cite{juhisz1991}  established the asymptotic behavior of  the
algebraic connectivity of ER random graph.
\begin{theorem}\cite{juhisz1991}\label{alg-ran}
Let $G$ be an ER random graph of $n$ vertices with the probability
$p$. Then for any $\varepsilon>0$, we have
$$ \alpha(G)=np+o(n^{1/2+\varepsilon}), {\rm in \  proabability}.
$$
\end{theorem}
 Further, Chung etc. (see \cite{chung 2003b} or \cite{chung 2006}) studied the semi-circle law for Laplacin
 eigenvalues of graphs.

The synchronization problem in networks of coupled oscillators is
closely related to  the consensus problem for network dynamics. The
consensus problems are related to the connections between spectral
properties of complex networks and ultrafast solution to distributed
decision-making problems for interacting groups of agents. Hence the
algebraic connectivity of networks is (locally) a measure of speed
of synchronization.

From mathematical view, the measurement of convergence speed of
solving the consensus problems in networks  is used by  the
algebraic connectivity of a network. There are some strongly
numerical evidences which support the conjecture that the network
dynamics on the small-world networks would display enhanced global
coordination compared to the regular lattices. In the recent paper
of Olfati-Saber \cite{olfati2008}, it is observed that the algebraic
connectivity of the small-world networks could be increased
dramatically by more than $1000$ times.  Gu etc. in \cite{gu2010}
gave a mathematical rigorous estimation of the lower bound for the
algebraic connectivity of the small-world networks, which is much
larger than the algebraic connectivity of the regular circle. This
result explains  why the consensus problems on the small-world
network have a ultrafast convergence rate and how much it can be
improved.

 A $2k$
regular lattice ${\mathcal{C}}(n, k)$ with $n$ vertices can be
constructed from a cycle of $n$ vertices and connected each vertex
to its $2k$ nearest neighbors. Let ${\mathcal{S}}(n, c, k)$ be a
small-world network that is  the union of a random graph
${\mathcal{G}}\big(n, \frac{c}{n}\big)$ and a $2k$ regular lattice
${\mathcal{C}}(n,k)$. Denote by $\lambda_2(n, c, k)$  the algebraic
connectivity of ${\mathcal{S}}(n, c, k)$. If $c = 0$, then
$\lambda_2(n, 0, k)$ is the algebraic connectivity of the $2k$
regular lattice ${\mathcal{C}}(n,k)$. The relationship between
parameter $c$ and the average shortcut per node $s$ is $c=2s$.

\begin{theorem}\cite{gu2010}
\label{main} Let ${\mathcal{S}}(n,c,k)$ be the small-world network
with $n$ nodes, which is a union of an Erd\"{o}s-R\'{e}yni random
graph ${\mathcal{G}}\big(n, \frac{c}{n} \big)$ and a $2k$ regular
cycle. Then the algebraic connectivity of ${\mathcal{S}}(n,c,k)$ is
almost surely bounded below by $$ \frac{k^2c^2\log\log
n}{2(k+1)^2\log^3 n}.
$$
\end{theorem}
Olfati and Saber \cite{olfati2008} defined
$$ \gamma_2(n, c, k) =
\frac{\lambda_2(n, c, k)}{\lambda_2(n, 0, k)}
$$
to be the {\it algebraic connectivity gain} of ${\mathcal{S}}(n, c,
k)$.

\begin{theorem}\cite{gu2010}
\label{gain} The algebraic connectivity gain of the small-world
network ${\mathcal{S}}(n,c,k)$ follows almost surely inequality
$$
\gamma_2({\mathcal{S}}(n,c,k)) \geq \frac{3kc^2 n^2 \log \log n}{2(k
+ 1)^3 (2 k+1) \pi^2 \log^3 n}. $$
\end{theorem}

\begin{center} \vskip 0.3cm {\bf Acknowledgements}
\end{center}
 Xiao-Dong Zhang would like to express their sincere gratitude to Professor Jiong-Sheng Li who
 introduced him the area of  spectral graph theory.

\end{document}